%% file: lanl27aug02.tex
\documentclass[11pt]{article}

\usepackage{epsfig,amsmath,latexsym}
\usepackage{amsfonts,amssymb,graphicx}
\def\slfrac#1#2{\hbox{\kern.1em %
  \raise.5ex\hbox{\the\scriptfont0 #1}\kern-.11em %
  /\kern-.15em\lower.25ex\hbox{\the\scriptfont0 #2}}}

\newtheorem{theorem}{Theorem}[section]
\newtheorem{lemma}[theorem]{Lemma}

\def\slfrac#1#2{\hbox{\kern.1em %
  \raise.5ex\hbox{\the\scriptfont0 #1}\kern-.11em %
  /\kern-.15em\lower.25ex\hbox{\the\scriptfont0 #2}}}

\newcommand{\beq}{\begin{eqnarray}}
\newcommand{\eeq}{\end{eqnarray}}
\newcommand{\beql}[1]{\begin{eqnarray}\label{#1}}
\newcommand{\beqs}{\begin{eqnarray*}}
\newcommand{\eeqs}{\end{eqnarray*}}
\newcommand{\eqn}[1]{(\ref{#1})}

\newcommand{\RR}{{\mathbb{R}}}
\newcommand{\ZZ}{{\mathbb{Z}}}
\newcommand{\sG}{{\mathcal{G}}}

\newcommand{\sS}{{\Sigma}}

\newcommand{\bv}{{\bf{v}}}
\newcommand{\bw}{{\bf{w}}}
\newcommand{\bx}{{\bf{x}}}
\newcommand{\by}{{\bf{y}}}

\newcommand{\qed}{\vrule height .9ex width .8ex depth -.1ex}

\makeatletter
\def\@sect#1#2#3#4#5#6[#7]#8{\ifnum #2>\c@secnumdepth
      \def\@svsec{}\else
      \refstepcounter{#1}\edef\@svsec{\csname the#1\endcsname.\hskip .75em }\fi
      \@tempskipa #5\relax
       \ifdim \@tempskipa>\z@
         \begingroup #6\relax
           \@hangfrom{\hskip #3\relax\@svsec}{\interlinepenalty \@M #8\par}%
         \endgroup
        \csname #1mark\endcsname{#7}\addcontentsline
          {toc}{#1}{\ifnum #2>\c@secnumdepth \else
                       \protect\numberline{\csname the#1\endcsname}\fi
                     #7}\else
         \def\@svsechd{#6\hskip #3\@svsec #8\csname #1mark\endcsname
                       {#7}\addcontentsline
                            {toc}{#1}{\ifnum #2>\c@secnumdepth \else
                              \protect\numberline{\csname the#1\endcsname}\fi
                        #7}}\fi
      \@xsect{#5}}
\def\@begintheorem#1#2{\it \trivlist \item[\hskip \labelsep{\bf #1\ #2.}]}

\def\plain{plain}\ifx\fmtname\plain\csname fi\endcsname
      
      \input docstrip
      \preamble

      Do not distribute the stripped version of this file.
      The checksum in the header refers to the documented version.

      \endpreamble
      \generateFile{here.sty}{t}{\from{here.doc}{}}
      \endinput
\fi
\ifcat a\noexpand @\let\next\relax\else\def\next{%
     \documentstyle[here,doc]{article}\MakePercentIgnore}\fi\next
\ifx\@Hxfloat\@Hundef\else\expandafter\endinput\fi
\let\@Hxfloat\@xfloat
\def\@xfloat#1[{\@ifnextchar{H}{\@HHfloat{#1}[}{\@Hxfloat{#1}[}}
\def\@HHfloat#1[H]{%
\expandafter\let\csname end#1\endcsname\end@Hfloat
\vskip\intextsep\vbox\bgroup\def\@captype{#1}\parindent\z@
\ignorespaces}
\def\end@Hfloat{\egroup\vskip \intextsep}

\makeatother

\catcode`\@=11

\catcode`@=12

\begin{document}

\begin{center}
{\large {\bf The Number of Triangles Needed to Span a Polygon Embedded in $\RR^d$}} \\
\vspace*{1.0\baselineskip}
{\em Joel Hass}\\
\vspace*{.2\baselineskip}
Department of Mathematics, \\
University of California, Davis \\
Davis, CA 95616 \\
\vspace*{1.5\baselineskip}

{\em Jeffrey C. Lagarias}\\
\vspace*{.2\baselineskip}
AT\&T Labs-Research, \\
Florham Park, NJ 07932-0971 \\
\vspace*{1.5\baselineskip}

(November 20, 2002)\\
\vspace*{1.5\baselineskip}
{\bf ABSTRACT}
\end{center}

Given a closed polygon $P$
having $n$ edges, embedded in $\RR^d$, we give upper
and lower bounds
for the minimal 
number of triangles $t$ needed to form a triangulated PL surface
embedded in $\RR^d$ having $P$ as its 
geometric boundary. More generally we
obtain such bounds for a triangulated 
(locally flat) PL surface having $P$ as its boundary
which is immersed in $\RR^d$ and whose interior is
disjoint from $P$. 
The most interesting case is dimension $3$, where the
polygon may be knotted. We use the
Seifert surface construction to show  that for any
polygon embedded in $\RR^3$ there exists
an embedded orientable triangulated PL surface
having at most $7n^2$ triangles, 
whose boundary is a subdivision of $P$. We complement
this with a construction of families of polygons
with $n$ vertices for which any such embedded surface requires at least
$\frac{1}{2}n^2 - O(n)$ triangles. We also
exhibit families of polygons in $\RR^3$
for which  $\Omega(n^2)$ triangles are required in any
immersed PL surface of the above kind.
In contrast, in dimension $2$ and in dimensions $d \ge 5$
there always exists an embedded locally flat  PL disk having
$P$ as boundary that contains at most $n$ triangles.
In dimension $4$ there 
always exists an immersed locally flat PL disk of
the  above kind that contains at most $3n$ triangles.
An unresolved case is that of embedded PL surfaces in dimension
$4$, where we establish only an  $O(n^2)$  upper bound.
These results can be viewed as providing 
qualitative discrete analogues of the isoperimetric 
inequality  for piecewise linear (PL) manifolds. In
dimension $3$ they
imply  that  the 
(asymptotic) discrete isoperimetric constant
lies between  1/2 and 7. \\


\section{Introduction}
\setcounter{equation}{0}

Given a closed polygon $P$ in $\RR^d$ 
having $n$ edges, we consider the
problem of giving  upper
and lower bounds
for the minimal 
number of triangles $t$ needed to form a triangulated PL surface
in $\RR^d$ having $P$ as its 
geometric boundary.
This has a well known  answer 
in $\RR^2$, which is 
$$
t = n - 2,
$$
in which $t$ is the number of triangles in any triangulation
of the (convex or nonconvex) polygon that adds no extra vertices,
see \cite[Theorem 23.2.1]{Su97}.
Such triangulations minimize the number of triangles 
over all possible triangulations, in which extra vertices might be added.

What happens in higher dimensions?
We consider two versions of the problem, in both of which
extra vertices are permitted, in the surface and added
to the polygonal boundary. 

(1) The surface is an
embedded oriented PL surface with no restriction on its genus.

(2) The surface is an immersed oriented
PL disk, with the extra restriction
that the interior of the surface cannot cut through its boundary.

\noindent Recall that a surface is {\em embedded} if it does not intersect itself, and
is {\em immersed} if it is locally embedded, i.e. if
each point on the surface has an embedded neighborhood.
More precisely, in  the
second case we allow only {\em complementary
immersed surfaces}, by which we mean immersed surfaces 
(of any genus) whose interior
does not intersect its boundary, and whose boundary is
embedded \footnote{This condition rules out
a surface whose boundary is a multiple covering of 
the given boundary curve.}. Complementary
immersed surfaces allow extra freedom over embedded surfaces, 
but in dimension $d=3$ they still detect
knottedness; a polygon $P$ in $\RR^3$ has 
a complementary immersed surface that is a topological
disk if and only if $P$ is unknotted. 
If we were to allow general immersed surfaces
rather than restricting to complementary immersed surfaces,
then the two-dimensional construction above works in all
higher dimensions, to produce an immersed disk
having $n-2$ triangles (whose interior in general will
intersect its boundary.) 
With the extra restrictions given above 
on the surfaces the answers becomes non-trivial.
We show that the minimal number of triangles grows like
$n^2$ in dimension $3$, is $O(n)$ in dimensions $5$ and
above, and also $O(n)$ in dimension $4$ for complementary immersed surfaces.

We first consider dimension $3$, where the
polygon may be knotted.
It is known that the Seifert surface construction
leads to an $O(n^2)$ algorithm  for constructing  oriented
embedded surfaces, see Vegter~\cite[p. 532]{Ve97}.
Using such a construction we  obtain the following explicit
upper bound for a triangulated oriented embedded surface
having the polygon as boundary,

\begin{theorem}~\label{upper-oriented}
For each closed polygonal curve
  $P$ with $n$ line segments embedded in $\RR^3$,
there exists an embedded oriented triangulated PL surface, having a
PL subdivision of $P$ as boundary, that has at most $7n^2$
triangles.
\end{theorem}

In this result, the polygonal curve $P$ is the geometric boundary
of the surface, but extra vertices may have to be added to
$P$ to get the boundary of the surface as a PL manifold.

Theorem~\ref {upper-oriented} shows that
for embedded surfaces of unspecified genus
the upper bound is polynomial in $n$. This result contrasts with bounds
for the size of an embedded oriented triangulated
PL surface of minimal genus spanning the curve $P$.
Hass, Snoeyink and Thurston~\cite{HST} show
in the unknotted case (minimal genus = $0$) that
$(C_1)^n$ triangles are sometimes required
as $n \to \infty$, where $C_1 >1$ is a fixed constant.
Complementing this result,
Hass, Lagarias and Thurston \cite{HLT} show that in the
unknotted case there always exists an embedded minimal
genus PL surface (an embedded disk) spanning the curve $P$,
that contains  at most $(C_2)^{n^2}$ triangles, for a fixed constant
$C_2 > 1$.

In \S3 we show that the upper bound of $O(n^2)$
in Theorem~\ref{upper-oriented} is the correct order
of magnitude. We present two constructions
based on different principles, giving $\Omega (n^2)$ lower bounds.

The first method involves the  genus $g(K)$ of the knot $K$.
The {\em genus} of a knot is the minimal genus of any orientable
embedded surface that has the knot as its boundary. The lower bound is
$$
t \geq 4g(K) + 1,
$$
and it applies to embedded orientable PL surfaces having  $K$
as boundary.
This lower bound depends only on the (ambient isotopy) type of the
knot $K$, so one gets the best result by minimizing the
number of edges in a polygon representing the knot, which
is called the {\em stick number} of the knot. Using
the $(n, n-1)$ torus knot,
we obtain the following result.

\begin{theorem}~\label{lower-oriented}
There exists an infinite sequence of values of $n \to \infty$
with  closed polygonal curves
  $P_n$ having  $n$ line segments embedded in $\RR^3$,
for which any
embedded triangulated PL surface, that is oriented and
that has a PL
subdivision of $P_n$ as boundary, requires at least
$\frac {n^2}{2} - 3n + 5$ triangles.
\end{theorem}

The second method uses an invariant of a knot diagram $K$,
the {\em writhe} $w(K)$. 
The writhe of a knot diagram $K$ is
obtained by assigning an orientation (direction) to the
knot diagram, and then assigning a sign of $\pm 1$ to each
crossing, with $+1$ assigned if the two directed paths of the
knot diagram at the crossing have the undercrossing oriented
by the right hand rule relative to the overcrossing, and
$-1$ if not. The writhe $w(K)$ of the oriented diagram
is the sum of these signs over all crossings; it is 
independent of the orientation chosen.

The lower
bound is
\begin{equation}~\label{wrbound}
t \geq  |w(K)|,
\end{equation}
and it applies to complementary immersed surfaces.
The writhe is not an ambient isotopy invariant, but is an invariant of a
of knot diagram
under Reidemeister moves of types II
and III only, with type I moves forbidden.

We apply this bound to show that there is
an infinite family of polygonal curves $P_n$ in $\RR^3$ having  a
quadratic lower bound for the number of triangles in
a  complementary immersed surface, of unrestricted genus. 

\begin{theorem}~\label{lower-oriented2}
There exists an infinite sequence of 
closed polygonal curves
  $P_n$ in $\RR^3$ having  $n$ line segments, with 
values of $n \to \infty$,
for which any
complementary immersed  triangulated PL surface
that has a PL
subdivision of $P_n$ as boundary, requires at least
$\frac {n^2}{36}$ triangles.
\end{theorem}

The writhe bound \eqn{wrbound} implies that a polygonal knot that
has a large writhe in one direction must have a large number
of crossings in any projection direction (Theorem~\ref{large.writhe}).

In \S4 we consider the
combinatorial isoperimetric
problem for embeddings of a curve in
dimensions $d \geq 4$. We obtain two
$O(n)$ upper bounds.
In these dimensions we construct  PL surfaces
spanning the polygon which are locally flat,
as defined at the beginning of \S4.
The local flatness condition is
a restriction on how the surface is situated in $\RR^d$.
It is known that local flatness always holds for an
embedded PL surface in codimension $3$ or more
(see \cite[Corollary 7.2]{RS82}), hence requiring local flatness
puts a constraint only in 
dimension $d=4$. 

We first treat dimension $d=4$, and construct a 
complementary immersed surface
that is locally flat and has $O(n)$ triangles.

\begin{theorem}~\label{four.dim}
Let $P$ be a closed polygonal curve
 embedded in $\RR^4$ 
consisting of $n$ line segments.
Then there exists a complementary immersed triangulated PL disk,
which is locally flat, has $P$ as its PL boundary, 
and contains  $3n$ triangles.
\end{theorem}

Second, in dimensions $d \ge 5$,
by coning the polygon $P$ to a suitable
point we obtain 
an  embedded PL disk with $n$ triangles. 

\begin{theorem}~\label{higher.dim}
Let $P$ be a closed polygonal curve
embedded in $\RR^d$,
with $d \ge 5$, consisting of $n$ line segments.
Then there exists an embedded triangulated PL disk
which is locally flat,  has
$P$ as its PL boundary, and  contains $n$ triangles.
\end{theorem}

To summarize, these results establish that the complexity of
the spanning surface is $O(n^2)$ in dimension $3$,
and is $O(n)$ in all other dimensions, except 
possibly in dimension $4$ for embedded surfaces.
The increased complexity in  dimension $3$
might be expected, since dimension $3$ is the only
dimension in which knotting is possible for curves.
As far as we know, the  remaining unresolved
case of embedded surfaces in
$\RR^4$ might conceivably have superlinear complexity;
if so, this would represent a new phenomenon
peculiar to the discrete case.
For this case we establish only an
$O(n^2)$ upper bound, as explained at the
end of \S5, while an $\Omega(n)$ lower bound is
immediate.

Our motivation for study of these questions comes from
an analogy with isoperimetric inequalities, which we considered 
in \cite{HST}.
The classical isoperimetric inequality asserts
that for a simple closed curve $\gamma$ of length $L$ in $\RR^2$, the
area $A$ that it encloses satisfies
$$
A \le \frac{1}{4\pi} L^2,
$$
with equality only in the case of a circle.
This inequality generalizes to all higher dimensions, where  
we allow either immersed surfaces, which can be restricted to be disks,
or embedded surfaces of arbitrary genus, as follows.
For a closed $C^2$-curve $\gamma$ of length $L$ embedded
in $\RR^d$ there exists an 
immersed disk of area $A$ having $\gamma$ as boundary,
as well as an embedded orientable surface of area $A$ having
$\gamma$ as boundary, such that in either case
$$
A \le \frac{1}{4\pi} L^2.
$$
The first of these $d$-dimensional results traces back to 
Beckenbach and Rado~\cite{BR33}, while the second
traces back to Blaschke~\cite{Bla30}, see
Osserman ~\cite[p. 1202]{Os76}. 
The problems we consider here are
discrete analogues of these two variants of the 
isoperimetric inequality.
The discrete measure of ``length''
of the polygon is the
number of line segments $n$ in its boundary, and the discrete measure
of the ``area'' of a
triangulated surface is the number of triangles  $t$ that it contains.
This type of combinatorial minimal area problem is associated
to affine geometry because 
these measures of ``length'' and ``area'' are both
affine invariant. It follows that our results are
most appropriately viewed
as results concerning
 $d$-dimensional affine space ${\mathbb A}^d$
without a metric structure, rather than
$\RR^d$ with its Euclidean structure.  However for convenience we 
formulate all results in $\RR^d$.

Our results determine the order of growth
of the discrete isoperimetric bounds as a function of $n$.
The discrete problem has some differences from the
classical problem, in that its bounds grow linearly
in $n$ rather than quadratically as in the classical isoperimetric
inequality, except in dimension $3$, and possibly
dimension $4$ for embedded surfaces. Our bounds 
are qualitative, so are not a
perfect analogue of the classical isoperimetric inequality
which gives an exact constant. For exact answers in the
discrete case there are an 
infinite number of cases, one for each value of the number
of edges $n$ in the polygon. It therefore
seems more natural to consider
a notion of asymptotic isoperimetric constant as $n \to \infty$,
We formulate this in the most interesting case of dimension $3$.
For each $n \ge 3$ we define the {\em discrete isoperimetric constant}
$\gamma (n)$ by
$$ 
\gamma(n) = \max_{P_n} \left( \min_{\Sigma ~\mbox{spans}~ P_n} 
~\frac{1}{n^2} t(\Sigma) \right),
$$
in which $P_n$ runs over all polygons with $n$ edges embedded in $\RR^3$,
and all surfaces $\Sigma$ are embedded surfaces.
We define the 
 {\em asymptotic discrete isoperimetric constant} $\Gamma$ 
in $\RR^3$ to be
$$ 
\gamma := \limsup_{n \to \infty}  \gamma (n). 
$$
Combining Theorems~\ref{upper-oriented} and
\ref{lower-oriented}
 implies that the asymptotic isoperimetric
constant $\tau$  must lie between 1/2 and 7.
It would be interesting  to determine 
whether the constant $\gamma$ is a limiting value
rather than a $\limsup$, and to determine its exact value.

A further direction for such  PL
isoperimetric problems
would be to establish isoperimetric bounds for
higher-dimensional submanifolds.
Consider a $k$-dimensional triangulated closed
PL-manifold $M$ embedded in $\RR^d$, where $k \geq 2$,
and ask: what is  the minimal number of $(k+1)$-simplices in an
embedded triangulated PL $(k+1)$-dimensional manifold
having a PL-subdivision of $M$ as its boundary?

Earlier work on  the complexity of embedded surfaces bounding 
unknotted curves in $\RR^3$ under
various restrictions includes Almgren and Thurston~\cite{AT77}.
Connections between combinatorial
complexity of such surfaces and the computational complexity 
of problems in knot theory appear in  \cite{HL}, \cite{HLP}.

%
%

\section{Upper Bound}
\setcounter{equation}{0}

We establish Theorem~~\ref{upper-oriented} by a straightforward analysis
of the construction due to Seifert~\cite{Se34} of
an orientable surface having a given knot as boundary.
A general description of Seifert surfaces and their construction
appears in  Rolfsen~\cite[Chapter 5]{Ro76}. \\

\noindent{\bf Proof of Theorem~\ref{upper-oriented}:}
Given a closed polygon $P$ in $\RR^3$ having $n$ line
segments, we first choose an orientation for it.
We obtain a knot diagram by orthogonally
projecting it onto a plane. Fix once and
for all  a projection direction
in ``general position'', so that  the projections of
any two line segments in $P$ intersect in at most one
point, and if the two segments in $P$ are disjoint
then this point must correspond to interior points of the
two segments. Without loss of generality we may rotate the
polygon so that the projection direction is in the $z$-direction
and the projected plane is $z= 0$, and we may translate it in
the $z$-direction so that it lies in the half-space $z \geq 1.$
The  projected image of the polygon in the plane
has $n$ vertices and $c$ crossing
points, where
$$
c \le n(n-3)/2,
$$
since an edge cannot intersect
its two adjacent edges or itself.
We make the projection
into a planar graph by marking vertices at each crossing point,
which we call {\em crossing vertices}.
This graph is a directed graph, with directed edges
obtained by projection of the orientation assigned to
the polygonal knot $P$, and is regarded as sitting in
the plane $z= 0$. Each vertex of this planar graph has either two
or four edges incident on it, so the faces of the graph can
be two-colored; call the
colors white and black.  The graph has
a single unbounded face, which we consider colored white;
it also has  at least one  bounded region colored black.
  We denote this
directed colored graph $\sG$; it has $n + c$ vertices, which
we call {\em initial vertices} in what follows.

We now add new vertices
to this graph as follows: At each edge containing
a crossing vertex we insert new vertices
very close to each of its crossing vertex
endpoints; call these {\em interior vertices}.
The resulting graph has $n + c$ {\em initial vertices} and
$4c$ interior vertices.  Each crossing vertex
has four edges incident to it.

Near each crossing vertex we now add edges
connecting pairs of interior vertices on adjacent
edges in cyclic order around the crossing point.
There are four such edges which
form a small quadrilateral enclosing the crossing vertex.
We choose the interior points close enough to the
crossing vertex so that the interior
of each edge in the boundary of this quadrilateral
does not intersect any other edge of the graph, and so
lies entirely inside one of the colored polygons of the original
graph. We assign that color to the edge. These four edges form two 
white edges and two black
edges; we discard the black edges and add the two white edges only, to form
an augmented planar graph.

The example of the
trefoil knot is pictured in Figure 1;
part (b) shows
the added vertices of the augmented planar graph, and the white
edges are indicated by dotted lines in (b).

The added white-colored
edges create two white-colored triangular regions adjacent
to each crossing vertex, which together form a ``bow-tie'' shaped region.
If one now deletes all crossing vertices and the four edges incident to
them (which have as other endpoint an interior vertex), and adds
in the white-colored edges only, then one obtains a new planar graph
$\sG'$, that has $n + 4c$ vertices.
This new graph may be disconnected, and consists of a union
of simple closed polygons. Its regions are
  two-colorable, with the coloring
obtained from that of $\sG$ by changing the color of the bow-tie
shaped regions from white to black.
For the trefoil knot this is pictured in Figure 1(c).

The graph $\sG'$
is a union of simple closed curves $C$ which we call {\em circuits},
some of which may be
nested inside others. To each circuit we assign an integer
level that measures its nesting.
An {\em innermost circuit} is one that contains
no other circuit: we assign  these innermost circuits level $1$. We
now inductively define a level for each other circuit, to be
one more than the maximal level of any circuit they contain;
the maximal level is at most $n$.

%
%

\clearpage
\begin{figure}[htb]
\begin{center}
\input fg31.pstex_t
\caption{Trefoil knot}
\end{center}
\label{trefoil.fig}
\end{figure}
\clearpage

We now construct a triangulated
embedded spanning surface $\Sigma$ for $\gamma$ as follows.

(1) For each circuit $C$  of level $k$ in $\sG'$ 
we make a copy $\tilde{C}$ of this
circuit in the plane $z = -k,$ i.e. we translate it
from the plane $z= 0$ by
the vector $(0, 0, -k)$. It forms a simple closed polygon
whose interior in this plane
 will form part of the surface $\sS$. If it has $m$ sides,
then we may triangulate it using $m - 2$
triangles lying in the plane $z= -k$.

(2) We next add vertical faces connecting the circuit $\tilde{C}$ to
the the polygon $P$ lying above it. More precisely,  we must
 enlarge $P$ to a one-dimensional simplicial complex $P'$
which includes preimages of the white edges.
We first add new vertices on $P$ which lie
vertically above the endpoints of the white edges;
these points are unique because the projection is one-to-one
off crossing points; this gives a subdivision of $P$.
We next  add
new edges (not on $P$) connecting these points, which
project to the white edges; $P'$ is the one-dimensional
simplicial complex in $\RR^3$ resulting from adding
these points and edges to $P$. To each
edge of the circuit $\tilde{C}$ there is a unique edge of
$P'$ that projects vertically onto it.
We take the convex hull of these two edges, that forms
a trapezoid with two vertical edges, plus the two edges
we started with.
These trapezoids will form part of the surface $\sS$.
Each of them may be triangulated by adding a diagonal
to the trapezoid. See Figure 2.

%
%

\begin{figure}[htb]
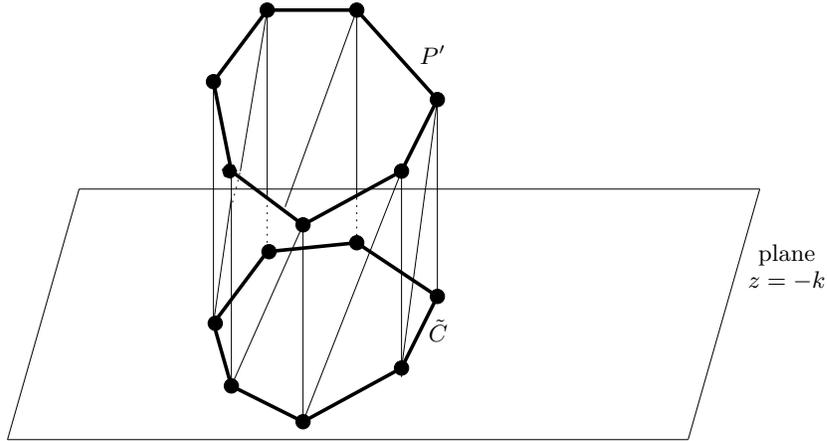

\begin{center}
\input fig2a.pstex_t
\end{center}
\begin{center}
\caption{Triangulated bowl-shaped region for cycle $C$}
\end{center}
\label{bowl-a}
\end{figure}

The part of the  surface $\sS$ associated to the
circuit $\tilde{C}$ in steps (1) and (2)  forms a bowl-shaped region whose
base is a  polygon in the plane $z= -k$
and part of $P'$ above $\tilde{C}$ as its lip.

(3) Above each bow-tie shaped region of $\sG'$ containing
a crossing vertex and two white edges, there lie four
edges of $P'$, two edges of which
project to the white edges on the plane $z=0$,
and the other two of which are part of edges of
$P$  whose projections on the
plane $z =0$ are disjoint except
at the crossing vertex. Let the vertices of the two
edges of $P'$ lying above the
white edges be labeled $[\bx_1, \bx_2]$ and $[\by_1, \by_2]$
respectively, with the black edges (not part of $P'$)
being $[\bx_1, \by_2]$
and $[\bx_2, \by_1]$,
and with the line segments $[\bx_1, \by_1]$
and  $(\bx_2, \by_2)$ being subsets of the original polygon $P$.
We then form the two triangles
$[\bx_1, \bx_2, \by_1]$ and $[\bx_2, \by_1, \by_2]$,
that share a common black edge $[\bx_2, \by_1]$, and add them to $\Sigma.$
See Figure 3.

%
%

\begin{figure}[htb]
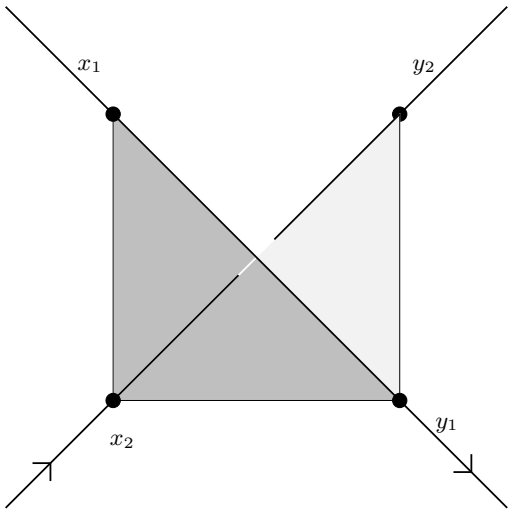

\begin{center}
\input fig2b.pstex_t
\end{center}
\begin{center}
\caption{Triangulated ``bow-tie'' region}
\end{center}
\label{bow-tie}
\end{figure}

We claim that $\Sigma$ forms a triangulated surface embedded in $\RR^3$,
which is orientable and has a
subdivision of $P$ as its boundary.

To see that $\sS$ is embedded in $\RR^3$, note that
the triangulated pieces (1)-(3) of  $\sS$  are embedded, and when projected
to the $z$-axis have disjoint interiors.
Thus these pieces can only overlap along
their boundary edges.

We use Seifert's argument to show $\sS$ is orientable. The contribution of
(1) and (2) corresponding to each circuit is a bowl-shaped surface
that is topologically a disk, with the crossing points located
near the lip of the bowl. At each crossing point the cup is
attached to another bowl by a rectangular strip with a half-twist
in it, twisting through an angle of $\pi$.  Since it is constructed
from disks attached along boundary intervals, $\sS$ is topologically
a $2$-manifold with boundary. In addition, the construction connects a
bowl at level $j$ only to bowls at level $j \pm 1$.
A compatible orientation then takes as
one side the upper (inside) surface of bowls at level $2j$ and
the lower (outside) surface of bowls at level $2j + 1$, plus
corresponding sides of the strips connecting them. Thus $\sS$ is orientable.

We bound the number of
triangles in $\sS$. The totality of triangles produced in
step (1) above is at most the number of
edges in all the circuits; this is at most the number
of edges in $P'$ after adding internal vertices, and is at most $n + 4c$.
In step (2) each trapezoid is associated to one of the at most
$n + 4c$ edges of the circuits, and has two triangles; thus these
contribute at most $2n + 8c$ triangles.
In step (3) there are two triangles added for each crossing vertex,
which totals $2c$. Thus the total is at most $3n + 14c$,
which is at most $7n^2 - 18n$ triangles.
This gives the desired upper bound $7n^2$. 
$~~~\qed$

\noindent {\bf Remark.} 
It is possible to modify the  construction in Theorem~\ref{upper-oriented}
to further improve the asymptotic upper bound to $cn^2$ with a constant
$c$ smaller than $7$. We leave the problem of obtaining the optimal
constant unanswered. 

%
%

\section{Lower Bound Constructions}
\setcounter{equation}{0}

We present two different constructions giving quadratic lower bounds
for triangulated surfaces.

\begin{lemma}\label{genus.bound}
Let $P$ be  a closed polygonal curve
embedded in $\RR^3$ whose associated
knot type $K$ has genus $g(K)$. If $t$ is
the number of triangles in a
triangulated oriented PL surface $\sS$ which is embedded
in $\RR^3$ and has a subdivision of $P$ as
boundary, then
\begin{equation}~\label{eq301}
t \geq 4g(K) + 1.
\end{equation}
\end{lemma}
{\bf Proof:}
Without loss of generality
we  may assume that the surface
$\sS$ is connected, by discarding any components
having empty boundary; this only decreases $t$.
If $V, E$ and $F$ denote the number of vertices, edge and faces
in the triangulated orientable surface $\sS$ , 
then its Euler characteristic is
$$\chi(\sS) = V - E + F.$$
By definition of knot genus this surface is of genus $g \ge g(K)$.
Recall that the genus of a surface with boundary $\sS$ is the smallest
genus of a connected surface $\sS'$ without boundary 
in which it can be embedded;
In this case $\sS'$ is obtained by gluing in a disk attached to
the (topological) boundary $P$. 
This adds one face, and no new edges or vertices,
hence we obtain
$$\chi(\sS) = -1 + \chi(\sS')= 1 - 2g \geq  1 - 2 g(K). $$
We have $F \geq t$, and since all faces in the surface are triangles,
we obtain
$$3t = 2E - m,$$
where $m$ is the number of edges on the boundary of the surface.
Counting the number of edges on the boundary gives
$$V \geq m.$$
 From these bounds follows
$$\chi(\sS) \geq  1 - 2g(K) = V - E + F \geq m - (\frac{3}{2}t + \frac{m}{2})
+ t,$$
which simplifies to
$$\frac{t}{2} \geq 2g(K) - 1 + \frac{m}{2} \geq 2g(K)+ \frac{1}{2},$$
since $m \geq 3$.
$~~~\qed$ \\

\noindent{\bf Remark.}
Define the
{\em unoriented genus} $g^*(K)$
of a knot $K$ to be the minimal value possible of
$1 - \frac{1}{2}\chi(\sS) $ taken over all embedded connected surfaces $\sS$,
orientable or not,  having
$K$ as boundary. Then $g^*(K)$ is an integer or half-integer, and
the same reasoning as above shows that 
$$ t \geq 4g^*(K) + 1$$
for the number of triangles in any triangulated PL surface
bounding a polygon $P$ of knot type $K$.

\noindent{\bf Proof of Theorem~\ref{lower-oriented}:} \\
We consider the $(m, m-1)$ torus knot $K_{m, m-1}$.
This has a polygonal representation
using $n=2m$ line segments, given in
Adams et al \cite[Lemma 8.1]{AB97}.  They also show that a
polygonal realization of this knot requires at least $2m$ segments
\cite[Theorem 8.2]{AB97}.

In 1934 Seifert~\cite[Satz 4]{Se34}
showed that the $(p, q)$-torus knot $K_{p,q}$ has genus
$$g(K_{p,q}) = \frac{(p-1)(q - 1)}{2}.$$
Thus we have $\displaystyle g(K_{m,m-1}) = \frac{m^2 - 3m + 2}{2}.$

We apply  Lemma~\ref{genus.bound} to $K_{m, m-1}$ and obtain
$t \ge 2m^2 - 6m + 5 \ge \frac{1}{2}n^2 - 3n + 5,$
as asserted.$~~~\qed$

We next obtain a lower bound in terms of the {\em writhe}  (or {\em 
Tait number})
of a knot diagram $K$ associated to $P$ by planar projection.
The writhe of a knot diagram $K$ is 
calculated by assigning an orientation (direction) to the
knot diagram, then associating a sign of $\pm 1$ to each
crossing, using $+1$  if the two directed paths of the
knot diagram at the crossing have the undercrossing oriented
by the right hand rule relative to the overcrossing, and
$-1$ if not. The writhe $w(K)$ of the oriented diagram
is the sum of these signs over all crossings.
The quantity $w(K)$ is independent of
the orientation, but depends on the direction of projection.

\begin{lemma}\label{writhe.bound}
Let $P$ be  a closed polygonal curve
embedded in $\RR^3$ that has an orthogonal  planar projection $K$
that has writhe $w(K)$.
If $t$ is the number of triangles in a complementary immersed
PL surface
in $\RR^3$ which is triangulated
and has a subdivision of $P$ as
boundary, then
\begin{equation}~\label{eq302}
t \geq |w(K)| + 1.
\end{equation}
\end{lemma}
{\bf Proof:}
The quantity $|w(K)|$ reflects
the amount of twisting between two different  longitudes
of the knot, the $z$-pushoff and the preferred  longitude.
The {\em preferred  longitude} is a longitude
on the knot defined by an embedded two-sided surface bounding
the knot $P$. The preferred  longitude is defined intrinsically as the
two primitive homology
classes $\pm [\tau]$ of a peripheral torus of the knot
that are annihilated on injecting
into $H^{1}(\RR^3 - P, \ZZ)$  (a peripheral torus
is the boundary of an
embedded regular neighborhood of the knot).
The $z$-pushoff is the curve on the peripheral torus directly above $K$
in the $z$ direction.

Any triangulated complementary immersed
surface $\sS$ having a subdivision of $P$ as
boundary necessarily defines a curve on the peripheral torus
in the class $\pm [\tau]$. The total twisting about $K$ of the boundary of a
smooth surface $\sS$ having $P$
as boundary, relative to the $z$-direction, is
given by $\pi w(K) $.
In the case of a PL surface this
total twisting can be computed by adding successive jumps in
a normal vector to $P$ pointing along the surface as one travels
along the (subdivided) curve $P$.
Since triangles are flat, twisting  occurs only at
the boundaries of two adjacent triangles, and can
be no more than $\pi$ at such a point. 
While it is possible that a triangle meets $P$ at
as many as three points, the total twisting contributed
by any triangle at all points where it meets $P$ is
at most $\pi$, the sum of its interior angles. It follows that there must be
more than $|w(K)|$ triangles in the surface that meet
the polygon $P$. $~~~\qed$ \\

\noindent{\bf Proof of Theorem~\ref{lower-oriented2}:} \\
There exists  a family of 
polygonal curves $P_m$ for  $m \geq 1$ having
$n = 6m+3$ segments and with writhe $w(P_m) = m(m+1)$.
The knot diagram for the polygon $P_3$ is pictured
in Figure 4
below. The construction for
general $m$ consists of adding more parallel strands to the pattern.

%
%

\begin{figure}[htb]
\centerline{\psfig{file=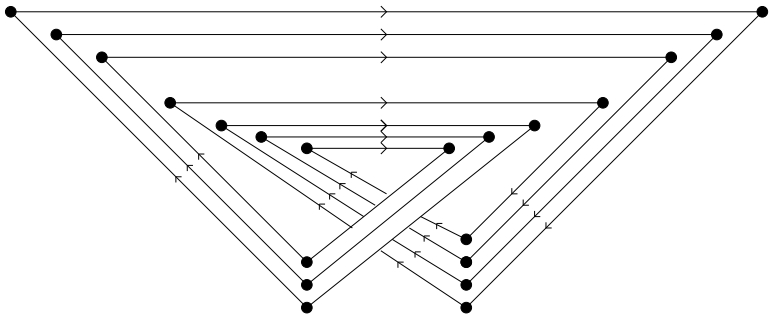,width=4in}}
\caption{Knot diagram  of polygon $P_3$.}
\label{bigwrithe}
\end{figure}

\noindent The theorem now
follows by applying the bound of Lemma~\ref{writhe.bound}, namely 
$\displaystyle t \ge \frac{n^2}{36} + \frac{3}{4}.$
$~~~\qed$ \\

Combining Lemma~\ref{writhe.bound} with the construction
of a surface in Theorem~\ref{upper-oriented}
yields a new result in knot theory. It says that if
a polygon $P$ embedded in $\RR^3$ has
a large writhe in some projection direction, relative to its
number of edges $n$, then in all projection
directions it has a large number of crossings.

\begin{theorem}\label{large.writhe}
Let $P$ be a polygonal knot embedded in $\RR^3$.
If $w(K)$ is the writhe
of one projection of $P$, then the number of crossings $c$ of any
projection $K'$ of $P$ satisfies
$$c  \ge \frac{1}{16}(|w(K)| - 3n). $$
\end{theorem}

{\bf Proof:}
By Lemma~\ref{writhe.bound} one has $t \ge |w(K)| + 1$,
However if a knot projection $K$ has $c$ crossings, then by the proof
of Theorem~\ref{upper-oriented} one can construct an oriented, embedded,
PL triangulated  surface having $P$ as
boundary with $t \le 3n + 14c$ triangles. Combining these
estimates gives the lemma.
$~~~\qed$

As an example, the polygons $P_m$ in $\RR^3$ 
given in the proof of Theorem~\ref{lower-oriented2}
(see  Figure~\ref{bigwrithe}) must
have crossing number
$ c \ge (m^2 - 17m - 9)/16$ in any projection.

\paragraph{Remarks.} 
(1) The  polygonal curves $P_m$ used in the proof of
Theorem~\ref{lower-oriented2}  are knotted.
We do not know how large $|w(P)|$ can be for an
unknotted polygon with $n$ edges. 
One can easily construct representatives of the unknot having writhe
$|w(P)| > cn$, for a positive constant $c$ and $n \to \infty$, but
we do not know whether it is possible to get representations $P$ of
the unknot having $n$ crossings and writhe $|w(P)| > cn^2$
with  $n \to \infty.$ 

(2) The proofs of Theorem~\ref{lower-oriented2} use topological
invariants associated to knottedness.  It may 
be that a quadratic lower bound can  hold for strictly geometric
reasons. 
A relevant geometric construction was given by
Chazelle ~\cite{Ch84}, who used it to
give examples of polyhedra with $n$ faces which
require $\Omega(n^2)$ tetrahedra in any triangulation.
A reviewer has suggested that this approach might
conceivably produce a family of unknotted polygons with
a $\Omega(n^2)$ lower bound.
Chazelle's construction  takes a hyperbolic
paraboloid, $H$ and makes two parallel translates $H^{+}$
and $H^{-}$ just above it and just below it. Now
$H$ is a ruled surface, and one takes $n$ line
segments connecting close points in  a ruling of $H$, and
$n$ segments each in the conjugate ruling of $H^{+}$ and $H^{-}$,
so that each pair of segments from opposite rulings cross
in vertical projection, so there are $\Omega(n^2)$ crossings
under vertical projections.
Then one connects the $3n$ segments in a zigzag manner to
produce a polygon with at most  $9n$ segments. 
(There is some freedom of choice in how to make the connections,
allowing the construction of polygons of various knot types.)
The geometric
principle to exploit is 
that the projection of any triangle with one edge on a segment
with endpoints in $H$ cannot cross more than a small number of projections
of segments in $H^{+}$ and $H^{-},$ unless the triangle is very narrow.
It seems plausible that an $\Omega(n^2)$  lower bound can be proved
for such a construction, but we leave this as an open problem.
This approach, if successfully
carried out,  would give lower bounds that apply to the complexity
of unoriented spanning surfaces.

%
%

\section{Higher Dimensions}
\setcounter{equation}{0}

In this section we consider polygons $P$ embedded in
$\RR^d$, for
dimensions $d \ge 4$, and construct locally flat PL surfaces having 
$P$ as boundary. Recall that a surface $\Sigma$ is 
{\em locally flat} if at
each point $\bx$ of the surface $\Sigma$ there is a neighborhood in 
$\RR^d$ homeomorphic
to $D \times I^{d - 2}$ in $\RR^d$, where $D$ is a topological $2$-disk in
the surface (or is a half $2$-disk with boundary for a boundary point $\bx$
of $\Sigma$)
and $I=[-1, 1]$ and $D \times {0}^{d-2}$ is part of $\Sigma$, 
see \cite[p. 36]{Ro76}, \cite[p. 50]{RS82}. 
For immersed surfaces we interpret local
flatness to apply to each sheet of the surface separately.

\noindent\paragraph{Proof of Theorem~\ref{four.dim}:}
We are given a closed polygon $P$ with $n$ edges embedded
in $\RR^4$, and vertices $\bv_1, ..., \bv_n$.
We can always pick a point $z \in \RR^4$
such that coning the polygon $P$ to the point $z$ will produce
a complementary immersed surface. However this surface need not
be locally flat at the cone point. To circumvent this problem,
we replace the cone point with a convex planar polygon $Q$ having $n$ vertices
and produce a  triangulated immersed
(polygonal) annulus connecting $P$ to $Q$.
Combining this with a triangulation of 
$Q$ to its centroid will yield
the desired locally flat immersed surface.

Given  a convex planar polygon $Q$ with vertices $\bw_1, ... , \bw_n$
we form the
\linebreak
(immersed) triangulated annulus $\Sigma_1$ between $P$ and 
$Q$ with triangles
\linebreak
$[\bv_{j}, \bv_{j+1}, \bw_{j+1}]$
and $[\bv_{j}, \bw_{j}, \bw_{j+1}]$, for
$1 \le j \le n $, using the convention that  $\bv_{n+1}:=\bv_1$
and $\bw_{n+1} :=\bw_1$. Then we triangulate $Q$ to its centroid
vertex $\bv_0$, obtaining a triangulated
disk $\Sigma_2$.  For ``general position'' $Q$ (described below)
this  construction produces an immersed
surface $\Sigma = \Sigma_1 \cup \Sigma_2$ 
consisting of  $n$ triangles from triangulating
$Q$ and $2n$ triangles in the annular part, for $3n$ triangles in all.
In general the surface $\Sigma_1$ is immersed rather than embedded,
because the triangles in it may intersect each other.

We show that $Q$ can be chosen so that $\Sigma$ is a complementary
immersed surface. The main problem is to ensure that $\Sigma_1$
intersects $P$ only in its boundary $\partial \Sigma_1$. We wish
$\Sigma_1$ to contain no line segment in any of its triangles
which intersects $P$ in two or more points. We define a ``bad set''
to avoid. Take the polygon $P$, extend its line segments to straight
lines $l_j$ in $\RR^4$, and call a point ``bad'' if it is on a line 
connecting any two points on the extended polygon. The ``bad'' set $B$
consists of a union of $n(n+1)/2$ sets $B_{ij}$,
in which $B_{ij}$ is the union of all lines
connecting a point  of  line $l_i$ to a point of line $l_j$,
with $1 \le i < j \le n$. Each $B_{ij}$ is either a hyperplane
(codimension $1$) in $\RR^4$ or is a plane (codimension $2$ flat) in $\RR^4$.
If a plane $F$ (codimension $2$ flat) is picked in ``general position''
in $\RR^4$ it will intersect $B$ in at most $n(n-1)/2$ lines and
points. In particular, such an $F$ contains 
a (two-dimensional)  open set $U$ not intersecting $B$
and disjoint from the convex hull of $P$. We
choose $Q$ to lie in this open set. Now 
$\Sigma_2$ is the convex hull of $Q$, which lies in $U$, so does
not intersect $P$. We claim 
that  $\Sigma_1$ intersects $P$ in $\partial \Sigma_1$.
Indeed any point $\bx$ in $\Sigma_1$ not on $P$ 
lies on  a line connecting a point of
$P$ to a point of $Q$, and this line contains at most one point
of $P$ because the point in $Q$ is not in the ``bad set'' $B$.
Since we already know of one point on $P$ on this line, which
is not $\bx$, the claim follows. We conclude that $\Sigma$ is
a complementary immersed surface.

We next show that $\Sigma$ is a locally flat (immersed)
surface.  We need only verify this at the vertices of
$\Sigma$. At the vertices $\bv$ of $P$ three triangles meet,
so locally the configuration is three-dimensional, and
local flatness holds in the three-dimensional
subspace around $\bw$ determined by the edges,
(as it does for any embedded polyhedral surface in $\RR^3$)
and this extends to local flatness in $\RR^4$ by taking
a product in the remaining direction. At a
vertex $\bw$ of $Q$ five triangles meet. 
However two of these
triangles lie in the plane $F$ of the polygon $Q$, hence
for determining local flatness we may disregard the edge 
into the interior of the polygon, and treat the vertex
as having four incident triangles. Suppose the remaining
$4$ edge directions leaving $\bw$ span a four dimensional space.
Take an invertible  linear transformation $L$ that maps these vectors
to $\bx_1=(1,1, 0, 0), \bx_2=(1, -1, 0, 0), \bx_3=(-1, -1,\epsilon,0),$ 
and $\bx_4= (-1, 1, 0, \epsilon),$ in cyclic order. Because each
angle in the convex polygon is less than $\pi$, we conclude
that the interiors of all four triangles project 
onto the positive linear combinations
of consecutive vectors, e.g. the first triangle
maps into the region $\lambda_1 \bx_1 + \lambda_2\bx_2$
with $\lambda_1, \lambda_2 \ge 0$.
Now projection on the first two coordinates
in this new coordinate system extends to a local homeomorphism
$U_1 \times I^2$ in a neighborhood of the vertex, and pulling
back by $L^{-1}$ gives the required local flat structure in a
neighborhood of the vertex. If instead the four edge directions
span a three-dimensional space, then the argument used for a
vertex of $P$ applies. 
Thus $\Sigma$ is locally flat at each  vertex of $Q$.
Finally the centroid vertex added to 
the polygon $Q$ is obviously locally flat, and
we conclude that $\Sigma$ is locally flat.

Finally we note that $\Sigma$ is a topological disk, since it
is two-sided and topologically is an annulus glued onto a disk.
$~~~\qed$ \\

\noindent\paragraph{ Proof of Theorem~\ref{higher.dim}:}
Given the polygon $P$ in $\RR^d$, for $d \ge 5$ we cone it to a suitably
chosen point $z \in \RR^d,$ chosen so that the coning is an
embedding. It suffices to choose a ``general position'' point,
because two planes (codimension $d-2$ flats)  in $\RR^d$ generically have
empty intersections.
The resulting surface $\Sigma$ has $n$ triangles, and is
a topological disk.

By a standard result, see  
Rourke and Sanderson \cite[Corollary 5.7 and ~Corollary 7.2]{RS82},
 this embedded surface is
locally flat. $~~~\qed$ \\

We conclude with some remarks concerning the unresolved
case of embedded surfaces in $\RR^4$
having a given polygon $P$ with $n$ edges as boundary. 
First, one can
always find such a triangulated surface using at most
$21n^2$ triangles, as follows. Take a projection
of the polygon $P$ into a hyperplane $H$, resulting
in a polygon $P^{*}$ in $H$, picking a projection  direction
such that the vertical surface $\Sigma_1$ connecting $P$ to $P^{*}$
is embedded.
Theorem~\ref{upper-oriented} gives a surface $\Sigma_2$
that uses at most $7n^2$ triangles which 
lies entirely in $H$, and has a subdivision of $P^{**}$ of $P^{*}$
as boundary and with $P^{**}$ having at most
$7n^2$ vertices. We obtain a
triangulated vertical surface $\Sigma_1$ connecting
$P$ to $P^{**}$ using at most $14n^2$ triangles, and
$\Sigma = \Sigma_1 \cup \Sigma_2$ is the required surface.
It  can be checked that this surface is locally flat.

Second, the immersed
surface constructed in Theorem~~\ref{four.dim} can
be converted to an embedded surface of higher genus
by cut-and paste, but we show that for some $P$ 
such a surface must contain
$\Omega(n^2)$ triangles. 
Recall that the {\em $4$-ball genus} of
a knot embedded in
a hyperplane in $\RR^4$ is the smallest genus
of any spanning surface of it
that lies strictly in a half-space of $\RR^4$ on one
side of this hyperplane.
If we start with
a polygon $P$ in $\RR^4$ that lies in a hyperplane, the
construction of Theorem~\ref{four.dim} will (in general) produce
an immersed surface lying in a half-space on one side of the
hyperplane, and a cut-and-paste construction will preserve
this property. (Note that cut-and-paste in $4$-dimensions
to replace two triangles intersecting in an interior point
with a non-intersecting set may result in eight triangles.)
As noted earlier, the  $(2n, 2n-1)$ torus knot has a polygonal
representation $P_n$ in a hyperplane using $4n$ line segments
(Adams et al \cite[Lemma 8.1]{AB97})
while a result of Shibuya~\cite{Sh86} implies
that its 4-ball genus is
at least $2n(2n-1)/8$. Applying
Lemma~\ref{genus.bound}, we conclude that the number of
triangles needed in an embedded orientable PL surface
of this type spanning
$P_n$ must grow quadratically
in $n$. This can be taken as (weak) evidence in support of
the possibility that for embedded surfaces in $\RR^4$
the best combinatorial isoperimetric bound may be
$O(n^2)$.

\section*{About Authors}

Joel Hass is at the Department of Mathematics, 
University of California, Davis, CA 95616;
\textsl{hass@math.ucdavis.edu}.

Jeffrey C. Lagarias is at AT\&T Labs-Research, 
180 Park Avenue,
Florham Park, NJ 07932-0971;
\textsl{jcl@research.att.com}.

\section*{Acknowledgments}
The authors are indebted to the reviewers for several 
useful remarks and suggestions, and for bringing some
references to our attention. 

Part of the work on this paper was done by the authors during
a visit to the Institute for Advanced Study. 
The first author was partially suppored by NSF grant DMS-0072348, 
and by a grant to the Institute of Advanced Study by AMIAS.

\end{document}

%% file: fg31.pstex_t
\begin{picture}(0,0)%
\includegraphics{fg31.pstex}%
\end{picture}%
\setlength{\unitlength}{2960sp}%
\begingroup\makeatletter\ifx\SetFigFont\undefined%
\gdef\SetFigFont#1#2#3#4#5{%
  \reset@font\fontsize{#1}{#2pt}%
  \fontfamily{#3}\fontseries{#4}\fontshape{#5}%
  \selectfont}%
\fi\endgroup%
\begin{picture}(7171,13405)(3233,-13099)
\put(5551,-3661){\makebox(0,0)[b]{\smash{\SetFigFont{11}{13.2}{\rmdefault}{\mddefault}{\updefault}\special{ps: gsave 0 0 0 setrgbcolor}(a)  Trefoil knot diagram ($n=7$, $c=3$)\special{ps: grestore}}}}
\put(5626,-8236){\makebox(0,0)[b]{\smash{\SetFigFont{11}{13.2}{\rmdefault}{\mddefault}{\updefault}\special{ps: gsave 0 0 0 setrgbcolor}(b)  Augmented graph\special{ps: grestore}}}}
\put(8851,-11611){\makebox(0,0)[lb]{\smash{\SetFigFont{11}{13.2}{\rmdefault}{\mddefault}{\updefault}\special{ps: gsave 0 0 0 setrgbcolor}$\mbox{level 1}$\special{ps: grestore}}}}
\put(5701,-13036){\makebox(0,0)[b]{\smash{\SetFigFont{11}{13.2}{\rmdefault}{\mddefault}{\updefault}\special{ps: gsave 0 0 0 setrgbcolor}(c)  Graph ${\mathcal G}'$\special{ps: grestore}}}}
\put(5476,-11011){\makebox(0,0)[lb]{\smash{\SetFigFont{11}{13.2}{\rmdefault}{\mddefault}{\updefault}\special{ps: gsave 0 0 0 setrgbcolor}level 0\special{ps: grestore}}}}
\end{picture}

%% file: fig2a.pstex_t
\begin{picture}(0,0)%
\includegraphics{fig2a.pstex}%
\end{picture}%
\setlength{\unitlength}{2960sp}%
\begingroup\makeatletter\ifx\SetFigFont\undefined%
\gdef\SetFigFont#1#2#3#4#5{%
  \reset@font\fontsize{#1}{#2pt}%
  \fontfamily{#3}\fontseries{#4}\fontshape{#5}%
  \selectfont}%
\fi\endgroup%
\begin{picture}(6537,3679)(2089,-3823)
\put(5626,-2986){\makebox(0,0)[lb]{\smash{\SetFigFont{9}{10.8}{\rmdefault}{\mddefault}{\updefault}{$\tilde{C}$}%
}}}
\put(5551,-661){\makebox(0,0)[lb]{\smash{\SetFigFont{9}{10.8}{\rmdefault}{\mddefault}{\updefault}{$P'$}%
}}}
\put(8626,-2311){\makebox(0,0)[b]{\smash{\SetFigFont{9}{10.8}{\rmdefault}{\mddefault}{\updefault}{plane}%
}}}
\put(8626,-2536){\makebox(0,0)[b]{\smash{\SetFigFont{9}{10.8}{\rmdefault}{\mddefault}{\updefault}{$z= -k$}%
}}}
\end{picture}

%% file: fig2b.pstex_t
\begin{picture}(0,0)%
\includegraphics{fig2b.pstex}%
\end{picture}%
\setlength{\unitlength}{2960sp}%
\begingroup\makeatletter\ifx\SetFigFont\undefined%
\gdef\SetFigFont#1#2#3#4#5{%
  \reset@font\fontsize{#1}{#2pt}%
  \fontfamily{#3}\fontseries{#4}\fontshape{#5}%
  \selectfont}%
\fi\endgroup%
\begin{picture}(4244,4244)(2454,-7058)
\put(6076,-3361){\makebox(0,0)[rb]{\smash{\SetFigFont{9}{10.8}{\rmdefault}{\mddefault}{\updefault}$y_2$}}}
\put(3451,-6511){\makebox(0,0)[b]{\smash{\SetFigFont{9}{10.8}{\rmdefault}{\mddefault}{\updefault}$x_2$}}}
\put(6076,-6361){\makebox(0,0)[lb]{\smash{\SetFigFont{9}{10.8}{\rmdefault}{\mddefault}{\updefault}$y_1$}}}
\put(3076,-3361){\makebox(0,0)[lb]{\smash{\SetFigFont{9}{10.8}{\rmdefault}{\mddefault}{\updefault}$x_1$}}}
\end{picture}